\documentclass[12pt,reqno]{amsart}
\usepackage{amsmath, amsfonts, amssymb, amsthm, hyperref}
\usepackage{bm}
\allowdisplaybreaks[4]
\def\l{\left}
\def\r{\right}
\def\bg{\bigg}
\def\({\bg(}
\def\){\bg)}

\def\f{\frac}

\def\eq{\equiv}

\def\Z{\mathbb Z}

\def\N{\mathbb N}

\def\1{{\bf 1}}

\theoremstyle{plain}
\newtheorem{theorem}{Theorem}[section]
\newtheorem{lemma}{Lemma}

\newtheorem{conjecture}{Conjecture}

\theoremstyle{definition}

\theoremstyle{remark}
\newtheorem{remark}{Remark}

\def\<{\langle}
\def\>{\rangle}

\begin{document}
\title{On a conjectural congruence of Guo}

\begin{abstract}
Let $p\equiv3\pmod{4}$ be a prime and $r$ a positive integer. We show that
$$
\prod_{k=1}^{(p^{2r}-1)/2}\frac{4k-1}{4k+1}\equiv1\pmod{p^2}.
$$
This confirms a recent conjecture of Guo \cite[Conjecture 1.4]{guo20}.
\end{abstract}

\author{Chen Wang}
\address{Department of Mathematics, Nanjing University, Nanjing 210093,
People's Republic of China}
\email{cwang@smail.nju.edu.cn}

\author{Hao Pan}
\address{School of Applied Mathematics, Nanjing University of Finance and Economics, Nanjing 210023, People's Republic of China}
\email{haopan79@zoho.com}

\maketitle

\section{Introduction}
\setcounter{lemma}{0} \setcounter{theorem}{0}
\setcounter{equation}{0} \setcounter{conjecture}{0}
In 2016, Swisher \cite[(H.3) with $r=2$]{swisher15} conjectured that
\begin{equation}\label{swisher}
\sum_{k=0}^{(p^2-1)/2}\f{(\f12)_k^3}{k!^3}\eq p^2\pmod{p^5}
\end{equation}
for any prime $p>3$ and $p\eq3\pmod4$.

Recently, Guo \cite[Theorem 1.2]{guo20} obtained two $q$-congruences related to \eqref{swisher} modulo $p^4$. For any positive integer $n\eq3\pmod4$, he proved that
\begin{gather}
\label{guo1}\sum_{k=0}^{(n^2-1)/2}\f{(1+q^{4k+1})(q^2;q^4)_k^3}{(1+q)(q^4;q^4)_k^3}q^k\eq\f{[n^2]_{q^2}(q^3;q^4)_{(n^2-1)/2}}{(q^5;q^4)_{(n^2-1)/2}}q^{(1-n^2)/2}\pmod{\Phi_n(q)^2\Phi_{n^2}(q)^2},\\
\label{guo2}\sum_{k=0}^{n^2-1}\f{(1+q^{4k+1})(q^2;q^4)_k^3}{(1+q)(q^4;q^4)_k^3}q^k\eq\f{[n^2]_{q^2}(q^3;q^4)_{(n^2-1)/2}}{(q^5;q^4)_{(n^2-1)/2}}q^{(1-n^2)/2}\pmod{\Phi_n(q)^2\Phi_{n^2}(q)^2}.
\end{gather}
Here $(a;q)_n$ is the Pochhammer symbol defined by $(a;q)_0$ and $(a;q)_n=(1-a)(1-aq)\cdots(1-aq^{n-1})$ for $n\geq1$, $[n]_q=1+q+\cdots+q^{n-1}$ denotes the $q$-integer, and
$$
\Phi_n(q)=\prod_{\substack{1\leq k\leq n\\ \gcd(n,k)=1}}(q-\zeta^k)
$$
is the $n$-th cyclotomic polynomial, where $\zeta$ is an $n$-th primitive root of unity.

Letting $n=p\eq3\pmod4$ be a prime and letting $q\to1$ in \eqref{guo1} and \eqref{guo2} we immediately obtain that
\begin{equation}\label{guo3}
\sum_{k=0}^{(p^2-1)/2}\f{(\f12)_k^3}{k!^3}\eq\sum_{k=0}^{p^2-1}\f{(\f12)_k^3}{k!^3}\eq p^2\prod_{k=1}^{(p^2-1)/2}\f{4k-1}{4k+1}\pmod{p^4}.
\end{equation}
If \eqref{swisher} holds, one must have
$$
\prod_{k=1}^{(p^2-1)/2}\f{4k-1}{4k+1}\eq1\pmod{p^2}.
$$
Motivated by this, Guo \cite[Conjecture 1.4]{guo20} proposed the following deeper conjecture.
\begin{conjecture}\label{guoconj}
Let $p\eq3\pmod4$ be a prime and $r$ a positive integer. Then
\begin{equation}\label{guo4}
\prod_{k=1}^{(p^{2r}-1)/2}\f{4k-1}{4k+1}\eq1\pmod{p^2}.
\end{equation}
\end{conjecture}

Our main purpose is to confirm the above conjecture.

\begin{theorem}\label{maintheorem1}
Conjecture \ref{guo4} is true.
\end{theorem}

\begin{remark}
Theorem \ref{maintheorem1} also indicates that \eqref{guo1} is actually a $q$-analogue of \eqref{swisher} modulo $p^4$.
\end{remark}

\section{Proof of Theorem \ref{maintheorem1}}
\setcounter{lemma}{0} \setcounter{theorem}{0}
\setcounter{equation}{0} \setcounter{conjecture}{0}
Recall that Morita's $p$-adic gamma function \cite{AAR,robert} is a $p$-adic analogue of the classical gamma function. For any integer $n\geq1$, the $p$-adic gamma function is defined by
$$
\Gamma_p(n):=(-1)^n\prod_{\substack{1\leq k<n\\ p\nmid k}}k.
$$
In particular, set $\Gamma_p(0):=1$. Let $\Z_p$ denote the ring of all $p$-adic integers. Since $\N=\{0,1,2,\ldots\}$ is dense in $\Z_p$ in the sense of $p$-adic norm $|\cdot|_p$, for any $x\in\Z_p$ define
$$
\Gamma_p(x):=\lim_{\substack{n\in\N\\ |x-n|_p\to0}}\Gamma_p(n).
$$
By the definition of $p$-adic gamma function over $\N$ and the continuity, we have
\begin{equation}\label{padicgamma1}
\f{\Gamma_p(x+1)}{\Gamma_p(x)}=\begin{cases}\displaystyle -x,\quad p\nmid x,\vspace{1.5mm}\\ \displaystyle -1,\quad p\mid x\end{cases}.
\end{equation}
It is known that for any $x\in\Z_p$,
\begin{equation}\label{padicgammaleg}
\Gamma_p(x)\Gamma_p(1-x)=(-1)^{p-\<-x\>_p},
\end{equation}
where $\<x\>_p$ denotes the least nonnegative residue of $x$ modulo $p$.
Clearly, \eqref{padicgammaleg} is a $p$-adic analogue of the well-known Legendre relation
$$
\Gamma(x)\Gamma(1-x)=\f{\pi}{\sin \pi x}.
$$

Two show Theorem \ref{maintheorem1} we also need the following known result which can be easily deduced from \cite[Theorem 14]{longrama}.

\begin{lemma}\label{long} For any prime $p\geq3,a,m\in\Z_p$ we have
$$
\Gamma_p(a+mp)\eq\Gamma_p(a)+\Gamma_p'(a)mp\pmod{p^2}.
$$
\end{lemma}

\medskip
\noindent {\it Proof of Theorem \ref{maintheorem1}}. We shall prove this result by induction on $r$. If $r=1$, then by \eqref{padicgamma1} we have
\begin{align*}
&\prod_{k=1}^{(p^{2}-1)/2}\f{4k-1}{4k+1}=\prod_{k=1}^{(p^{2}-1)/2}\f{k-1/4}{k+1/4}=\f{\l(\f34\r)_{\f{p^2-1}2}}{\l(\f54\r)_{\f{p^2-1}2}}=\f{\Gamma\l(\f{2p^2+1}{4}\r)\Gamma\l(\f54\r)}{\Gamma\l(\f34\r)\Gamma\l(\f{2p^2+3}{4}\r)}\\
=&\f{\f{p}{4}\cdot\f{5p}{4}\cdots\f{2p^2-p}{4}}{\f{3p}{4}\cdot\f{7p}{4}\cdots\f{2p^2-3p}{4}}\cdot\f{\Gamma_p\l(\f{2p^2+1}{4}\r)\Gamma_p\l(\f54\r)}{\Gamma_p\l(\f34\r)\Gamma_p\l(\f{2p^2+3}{4}\r)}=\f{p\l(\f14\r)_{\f{p+1}{2}}}{\l(\f34\r)_{\f{p-1}{2}}}\cdot\f{\Gamma_p\l(\f{2p^2+1}{4}\r)\Gamma_p\l(\f54\r)}{\Gamma_p\l(\f34\r)\Gamma_p\l(\f{2p^2+3}{4}\r)}\\
\eq&\f{p(-1)^{(p+1)/2}\Gamma_p\l(\f{2p+3}{4}\r)\Gamma_p\l(\f34\r)}{\f{p}{4}(-1)^{(p-1)/2}\Gamma_p\l(\f14\r)\Gamma_p\l(\f{2p+1}{4}\r)}\cdot\f{\Gamma_p\l(\f14\r)\Gamma_p\l(\f54\r)}{\Gamma_p\l(\f34\r)^2}\\
=&\f{\Gamma_p\l(\f{2p+3}{4}\r)\Gamma_p\l(\f14\r)}{\Gamma_p\l(\f{2p+1}{4}\r)\Gamma_p\l(\f34\r)}\pmod{p^2}.
\end{align*}
In view of \eqref{padicgammaleg} and Lemma \ref{long} we have
\begin{align*}
\f{\Gamma_p\l(\f{2p+3}{4}\r)\Gamma_p\l(\f14\r)}{\Gamma_p\l(\f{2p+1}{4}\r)\Gamma_p\l(\f34\r)}=&\Gamma_p\l(\f{2p+3}{4}\r)\Gamma_p\l(\f{3-2p}{4}\r)(-1)^{\f{p+1}4}\Gamma_p\l(\f14\r)^2(-1)^{\f{p+1}{4}}\\
\eq&\Gamma_p\l(\f34\r)^2\Gamma_p\l(\f14\r)^2=1\pmod{p^2}.
\end{align*}
Now we suppose that the result holds for $r<n$, then for $r=n$ we also have
\begin{align*}
&\prod_{k=1}^{(p^{2n}-1)/2}\f{4k-1}{4k+1}=\f{\l(\f34\r)_{\f{p^{2n}-1}2}}{\l(\f54\r)_{\f{p^{2n}-1}2}}=\f{\Gamma\l(\f{2p^{2n}+1}{4}\r)\Gamma\l(\f54\r)}{\Gamma\l(\f34\r)\Gamma\l(\f{2p^{2n}+3}{4}\r)}\\
=&\f{\f{p}{4}\cdot\f{5p}{4}\cdots\f{2p^{2n}-p}{4}}{\f{3p}{4}\cdot\f{7p}{4}\cdots\f{2p^{2n}-3p}{4}}\cdot\f{\Gamma_p\l(\f{2p^{2n}+1}{4}\r)\Gamma_p\l(\f54\r)}{\Gamma_p\l(\f34\r)\Gamma_p\l(\f{2p^{2n}+3}{4}\r)}=\f{p\l(\f14\r)_{\f{p^{2n-1}+1}{2}}}{\l(\f34\r)_{\f{p^{2n-1}-1}{2}}}\cdot\f{\Gamma_p\l(\f{2p^{2n}+1}{4}\r)\Gamma_p\l(\f54\r)}{\Gamma_p\l(\f34\r)\Gamma_p\l(\f{2p^{2n}+3}{4}\r)}\\
=&-\f{p\cdot\f{3p}{4}\cdot\f{7p}{4}\cdots\f{2p^{2n-1}-3p}{4}}{\f{p}{4}\cdot\f{5p}{4}\cdots\f{2p^{2n-1}-p}{4}}\cdot\f{\Gamma_p\l(\f{2p^{2n-1}+3}{4}\r)\Gamma_p\l(\f34\r)}{\Gamma_p\l(\f14\r)\Gamma_p\l(\f{2p^{2n-1}+1}{4}\r)}\cdot\f{\Gamma_p\l(\f{2p^{2n}+1}{4}\r)\Gamma_p\l(\f54\r)}{\Gamma_p\l(\f34\r)\Gamma_p\l(\f{2p^{2n}+3}{4}\r)}\\
=&\f{\l(\f34\r)_{\f{p^{2n-2}-1}{2}}}{\l(\f54\r)_{\f{p^{2n-2}-1}{2}}}\cdot\f{\Gamma_p\l(\f{2p^{2n-1}+3}{4}\r)}{\Gamma_p\l(\f{2p^{2n-1}+1}{4}\r)}\cdot\f{\Gamma_p\l(\f{2p^{2n}+1}{4}\r)}{\Gamma_p\l(\f{2p^{2n}+3}{4}\r)}.
\end{align*}
By the induction hypothesis, we have
$$
\f{\l(\f34\r)_{\f{p^{2n-2}-1}{2}}}{\l(\f54\r)_{\f{p^{2n-2}-1}{2}}}\eq1\pmod{p^2}.
$$
Therefore,
\begin{align*}
&\prod_{k=1}^{(p^{2n}-1)/2}\f{4k-1}{4k+1}\eq\f{\Gamma_p\l(\f{2p^{2n-1}+3}{4}\r)}{\Gamma_p\l(\f{2p^{2n-1}+1}{4}\r)}\cdot\f{\Gamma_p\l(\f14\r)}{\Gamma_p\l(\f34\r)}\eq1\pmod{p^2},
\end{align*}
where the last step follows from \eqref{padicgammaleg} and Lemma \ref{long} again.\qed

\end{document}